\documentclass[11pt,a4paper]{article}
\usepackage[utf8]{inputenc}
\usepackage[english]{babel}
\usepackage{amsmath,amssymb,amsthm}
\usepackage{geometry}

\numberwithin{equation}{section}

\newtheorem{theorem}{Theorem}[section]
\newtheorem{lemma}[theorem]{Lemma}
\theoremstyle{definition}
\newtheorem{remark}[theorem]{Remark}

\title{\textbf{Lehto-type growth conditions and a Picard theorem for the gradient of entire solutions to the Euler -- Lagrange equations of variational integrals with linear growth in the plane including a small extension of Moser's Bernstein result}}
\author{\textbf{Martin Fuchs}}
\date{}
\begin{document}
\maketitle
\begin{abstract}
We mainly consider entire solutions $u: \mathbb{R}^2 \rightarrow \mathbb{R}$ of the Euler -- Lagrange equation for a variational integral with density $g(|\nabla u|)$ of linear growth. Assuming a Lehto-type condition at infinity, we show that either $\nabla u$ is a constant vector or the image of $\nabla u$ is the whole plane with the possible exception of at most one point. This result is extended to certain energy densities of nearly linear growth. Moreover, we slightly generalize the Bernstein result of Moser.
\end{abstract}


\section{Introduction}

In our note we first discuss entire solutions $u \in C^2(\mathbb{R}^2)$ of the equation
\begin{equation}
\operatorname{div} \left( \frac{g'(|\nabla u|)}{|\nabla u|} \nabla u \right) = 0
\end{equation}
arising as the Euler equation of the variational integral $\int_{\Omega} g(|\nabla u|) \, dxdy$, where $\Omega$ is a domain in\\\ $\mathbb{R}^2$. Concerning the function $g: [0,\infty) \rightarrow \mathbb{R}$ we suppose that $g$ is of class $C^2([0,\infty))$ satisfying
\begin{equation}
g'(0) = 0, \quad g''(t) > 0,
\end{equation}
\begin{equation}
m(t-1) \le g(t) \le M(t+1)
\end{equation}
for all $t \ge 0$ with numbers $m, M > 0$. We further assume ($c$ denoting a positive constant)
\begin{equation}
\exists\mu > 1: g''(t) \ge c(1+t)^{-\mu} \quad \text{for all } t \ge 0
\end{equation}
and finally add the technical condition
\begin{equation}
t g''(t) \le g'(t), \quad t \ge 0.
\end{equation}
Evidently the second part of (1.2) is just a consequence of inequality (1.4), and (1.5) will follow, if for instance we require \[ g^{(3)}(t) \leq 0 \quad \text{for all } t \geq 0. \]\\

For the choice $g(t) = \sqrt{1+t^2}$ we have the nonparametric minimal surface case for which (1.4) holds with $\mu=3$, and as it was shown by Bernstein in his famous papers [1], [2], the only entire solutions are the affine functions $\mathbb{R}^2 \rightarrow \mathbb{R}$. This Bernstein property of equation (1.1) under the assumptions (1.2), (1.3) combined with the reverse inequality in (1.4) was extended independently by Mikljukov [3] and by Farina, Sciunzi and Valdinoci [4]: in [3] the author considers the case $\mu \ge 2$ actually working with the weaker condition
\[
\int_{0}^{\infty} t g''(t) \, dt < \infty
\]
replacing the opposite inequality in (1.4), whereas in [4] the authors can cover the case of exponents $\mu \ge 3$. The reader being in addition interested in the history of Bernstein's problem is also referred to the papers [3], [4], and in our note [5] we summarized both approaches in the setting of equation (1.1), but -- as we emphasize again -- with (1.4) being replaced by the opposite inequality.

To our knowledge the validity of Bernstein's theorem in the presence of condition (1.4) is an unsolved problem to which we can not present any substantial contribution. Quite the contrary, we start with the discussion of a rather strange behaviour of entire solutions $u$, whose gradients do not grow too fast at infinity: in this case either $u$ is an affine function or the image of $\nabla u$ is the whole plane with at most one exceptional point. More precisely, it holds

\begin{theorem} \label{theo1.1}
Let the conditions (1.2) -- (1.5) hold for the function $g$, where in (1.4) we recall that $\mu > 1$. Assume that $u \in C^2(\mathbb{R}^2)$ is an entire solution of (1.1) such that
\begin{equation}
|\nabla u(z)| \le c \left[ \ln (2+|z|) \right]^{\frac{1}{\mu-1}}, \quad z \in \mathbb{R}^2,
\end{equation}
with some constant $c>0$. Then we have
\begin{align*}
&\text{either : } \nabla u \equiv \text{const } \text{(\textit{``Bernstein case''})} \\
&\text{or : } \nabla u(\mathbb{R}^2) \supset \mathbb{C} - \{p\} \text{ for a point } p \in \mathbb{R}^2.
\end{align*}
\end{theorem}

\begin{remark}
In the second case of Theorem \ref{theo1.1} the situation is comparable to the so-called little Picard theorem from complex analysis (see, e.g.[6], Theorem 16.22, p. 357).
\end{remark}

\begin{remark}
We have no idea, if the ``or'' case in Theorem \ref{theo1.1} actually occurs. However,under rather strong assumptions, it can be excluded, we refer to Theorem 1.8 and Theorem 1.10 below.
\end{remark}

\begin{remark}
As we will discuss later, condition (1.6) can be replaced by the weaker requirement
\begin{equation}
\int_{1}^{\infty} \frac{1}{r K(r)} \, dr = \infty,
\end{equation}
where $K(r) := \int_{0}^{2\pi} \left( 1 + |\nabla u(r e^{i\theta})| \right)^{\mu-1} \, d\theta$. The divergence of the integral guarantees the validity of the famous Lehto-condition for quasiconformal mappings (see,e.g., [7]
or [8]) being of crucial importance for the proof of Theorem 1.1.\end{remark}

\begin{remark}
We note that in the situation of (1.6) the solution can not grow much faster than an affine function.
\end{remark}

Let us now briefly look at a more general setting: instead of using the radial Euler--Lagrange equation (1.1), we consider the more general divergence equation
\begin{equation}\label{eq:1}
\text{div}\,[Df(\nabla u)] = 0 \quad \text{on } \mathbb{R}^{2}
\end{equation}
with $u \in C^{2}(\mathbb{R}^{2})$ and an integrand $f \in C^{2}(\mathbb{R}^{2})$, $f = f(P)$. For some $\mu > 1$ and with constants $c_{1}, c_{2} > 0$, let the following ellipticity condition hold:
\begin{equation}\label{eq:2}
c_{1}(1+|P|)^{-\mu}|Q|^{2} \le D^{2}f(P)(Q,Q) \le c_{2}\frac{|Q|^{2}}{1+|P|}
\end{equation}
for all $P, Q \in \mathbb{R}^{2}$. Inequality (1.9) is just the condition of  $\mu$-ellipticity introduced in the paper [9]. We have in extension of Theorem 1.1 
\begin{theorem} \label{theo1.6}Let $u \in C^{2}(\mathbb{R}^{2})$ denote an entire solution of equation (1.8) with f satisfying (1.9). If in addition we have inequality (1.6), then the conclusions of Theorem 1.1. hold.\end{theorem}

\begin{remark} Clearly the conditions (1.2), (1.3),(1.4) and (1.5) imposed on the function g yield inequality (1.9) for the density $f(P):= g(t)$, where t is the norm of $P$.  We note that (1.9) inherits more flexibility, for instance, f not necessarily is of linear growth, however, the case of integrands with radial structure seems to be of more geometric interest. \end{remark}

Coming back to Bernstein's problem, we recall the famous result of Moser [10] stating that if in Theorem 1.1 and Theorem 1.6 condition (1.6) is replaced by the global boundedness of the gradient of the solution u, then u is an affine function. We have the following minor extension of Moser's theorem: 
\begin{theorem} \label{theo1.8}Let $u \in C^{2}(\mathbb{R}^{2})$ denote an entire solution of equation (1.8) with f satisfying (1.9). If inequality (1.6) is replaced by the condition \begin{equation}
    |\nabla u(z)| \leq c \ln(\ln|z|)
\end{equation}
with some constant $c > 0$ and for all points $z \in \mathbb{R}^2$ such that $|z|$ is sufficiently large, then $u$ is an affine function.\end{theorem}
\begin{remark} Condition (1.10) is much more restrictive than the upper bound from (1.6), however, from the proof of Theorem 1.8 we get \end{remark}
\begin{theorem}\label{theo1.10}Let the Landau $o$-condition for the gradient at infinity hold, i.e.
\begin{equation}
|\nabla u(z)| = o\left([\ln(2 + |z|)]^{\frac{1}{\mu-1}}\right) \quad \text{as } |z| \to \infty,
\end{equation}where $u \in C^{2}(\mathbb{R}^{2})$ denotes an entire solution of equation (1.8) with f satisfying (1.9). Then $u$ must be affine.\end{theorem}

Our paper is organized as follows: in Section 2 we transform (1.1) into a generalized Beltrami equation and collect some technical tools. The proof of Theorem \ref{theo1.1} is outlined in Section 3 making essential use of results due to Lehto and Virtanen presented in the monograph [7], and our line of arguments will mainly follow ideas taken from the work [8] of Astala, Iwaniec and Martin. In Section 4 we will add some short comments concerning the modifications necessary for the proof of Theorem 1.6, and Section 5 contains the proof of Theorem 1.8,  the proof of Theorem 1.10 being presented in Section 6.

\section{Transformation of equation (1.1) into a generalized Beltrami equation}

Consider $u \in C^2(\mathbb{R}^2)$ satisfying (1.1) and let the function $g$ just denote an element of $C^2([0,\infty))$ satisfying the condition (1.2). We define the quantities
\[
p := \frac{\partial}{\partial x} u = u_x, \quad q := u_y, \quad t := |\nabla u| = \sqrt{p^2+q^2},
\]
where $z = (x,y) = x+iy \in \mathbb{C}$, $u = u(x,y)$. Clearly we have the identity
\begin{equation}
p_y - q_x = 0,
\end{equation}
and (1.1) turns into
\begin{equation}
\partial_x \left( \frac{g'(t)}{t} p \right) + \partial_y \left( \frac{g'(t)}{t} q \right) = 0.
\end{equation}
From (2.1), (2.2) we infer
\begin{equation}
a_{11} p_x + 2 a_{12} p_y + a_{22} q_y = 0
\end{equation}
with coefficients defined through the formulas
\begin{equation}
\begin{cases}
a_{11} := \frac{1}{t} g'(t) + \frac{1}{t^2} p^2 \left[ g''(t) - \frac{1}{t} g'(t) \right], \\
a_{22} := \frac{1}{t} g'(t) + \frac{1}{t^2} q^2 \left[ g''(t) - \frac{1}{t} g'(t) \right], \\
a_{12} := a_{21} := \frac{1}{t^2} p q \left[ g''(t) - \frac{1}{t} g'(t) \right].
\end{cases}
\end{equation}

Quoting (1.2) it follows from (2.4) that
\begin{equation}
a_{11} + a_{22} = g''(t) + \frac{1}{t} g'(t) > 0
\end{equation}
holds, where in case $t=0$ the quantity $\frac{1}{t} g'(t)$ is just $g''(0)$.\\\\

Let us consider the function
\begin{equation}
\Phi(z) := \Phi(x,y) := p(x,y) - i q(x,y)
\end{equation}
and introduce the Wirtinger derivatives
\[
\partial_z := \frac{1}{2}(\partial_x - i\partial_y), \quad \partial_{\bar{z}} := \frac{1}{2}(\partial_x + i\partial_y).
\]

Then we have

\begin{lemma}
The function $\Phi$ from (2.6) is a solution of the generalized Beltrami equation
\begin{equation}
\partial_{\bar{z}} \Phi = K_1 \partial_z \Phi + K_2 \overline{\partial_z \Phi} \quad \text{on } \mathbb{C}
\end{equation}
with coefficients (recall (2.5))
\begin{equation}
K_1 = \frac{\frac{1}{2}(a_{22}-a_{11}) - i a_{12}}{a_{11}+a_{22}}, \quad K_2 = \frac{\frac{1}{2}(a_{22}-a_{11}) + i a_{12}}{a_{11}+a_{22}}.
\end{equation}
\end{lemma}

\begin{proof}[Proof of Lemma 2.1]
We express the quantities $p_x, q_y, p_y$ in terms of $\partial_z \Phi, \partial_{\bar{z}} \Phi, \overline{\partial_z \Phi}$ and return to (2.3) with the result
\[
\begin{aligned}
a_{11} & \left( \partial_{\bar{z}} \Phi + \frac{1}{2} \partial_z \Phi + \frac{1}{2} \overline{\partial_z \Phi} \right) + \\
2 a_{12} & \left( \frac{i}{2} \partial_z \Phi - \frac{i}{2} \overline{\partial_z \Phi} \right) + \\
a_{22} & \left( \partial_{\bar{z}} \Phi - \frac{1}{2} \partial_z \Phi - \frac{1}{2} \overline{\partial_z \Phi} \right) = 0.
\end{aligned}\]
After collecting terms this immediately yields equation (2.7) with $K_1$ and $K_2$ according to (2.8).
\end{proof}\
For preparing the next auxiliary result we remark that the eigenvalues of the symmetric matrix $(a_{ij})_{1 \le i,j \le 2}$ are given by
\begin{equation}
\lambda = g''(t) > 0, \quad \Lambda = g'(t)/t > 0
\end{equation}
being strictly positive on account of (1.2), hence the matrix of the coefficients $a_{ij}$ is positive definite.

\begin{lemma}
The determinant of the Jacobian matrix of $\Phi$ from (2.6) is nonnegative, i.e.
\[
\det D\Phi = \det \begin{pmatrix} p_x & -q_x \\ p_y & -q_y \end{pmatrix} \ge 0.
\]
\end{lemma}

\begin{proof}[Proof of Lemma 2.2]
It holds on account of (2.1)
\[
\det D\Phi = -p_x q_y + p_y q_x = p_y^2 - p_x q_y.
\]
Recalling (2.4) we have
\[
a_{22} = \frac{1}{t^2} q^2 g''(t) + \frac{1}{t} g'(t) \left( 1 - \frac{q^2}{t^2} \right) > 0,
\]
and (2.3) implies
\[
q_y = -\frac{a_{11}}{a_{22}} p_x - 2 \frac{a_{12}}{a_{22}} p_y,
\]
hence
\begin{equation*}
\begin{split}
\det D\Phi &= p_y^2 + \frac{a_{11}}{a_{22}} p_x^2 + 2 \frac{a_{12}}{a_{22}} p_x p_y \\
&= \frac{1}{a_{22}} \left( a_{22} p_y^2 + a_{11} p_x^2 + 2 a_{12} p_x p_y \right) \ge 0
\end{split}
\end{equation*}
which proves our claim.
\end{proof}
We finish this section by adding a few technical comments:
returning to the defining equations (2.8) it is immediate that $K_2 = \overline{K}_1$, moreover, we have the identity
\[
\begin{aligned}
|K_1|^2 & = \frac{\frac{1}{4}(a_{22} - a_{11})^2 + a_{12}^2}{(a_{11} + a_{22})^2} = \\
& \frac{\frac{1}{4}(a_{11} + a_{22})^2 - a_{11} a_{22} + a_{12}^2}{(a_{11} + a_{22})^2} = \\
& \frac{1}{4} - \frac{a_{11} a_{22} - a_{12}^2}{(a_{11} + a_{22})^2}
\end{aligned}
\]
hence
\[
|K_1|^2 = \frac{1}{4} - \frac{\det(a_{ij})}{(\operatorname{trace}(a_{ij}))^2} = \frac{1}{4} - \frac{\lambda \Lambda}{(\lambda+\Lambda)^2} < \frac{1}{4}
\]
on account of (2.9), and we find
\begin{equation}
|K_1| < \frac{1}{2}.
\end{equation}
\section{Proof of Theorem 1.1}
Let the conditions (1.2) -- (1.5) hold for the function $g$ and consider an entire solution $u \in C^2(\mathbb{R}^2)$ of equation (1.1) for which we have (1.6). We consider two cases: \\\\Case 1: $\Phi$ is a constant function.\\Then our claim is a direct consequence of equation (2.6).\\\\Case 2: $\Phi$ is not constant.\\In order to prove the Picard-property of the gradient mapping, we will apply techniques outlined in [8], Chapter 20.9, demonstrating the strength of the Lehto-condition . According to Lemma 2.1, the complex gradient mapping $\Phi(z) = p(x, y) - iq(x, y)$ is of class $C^1(\mathbb{R}^2)$ and satisfies the generalized Beltrami-type system (recall (2.7))
\begin{equation*}
\partial_{\bar{z}}\Phi(z) = K_1(z)\partial_z\Phi(z) + K_2(z)\overline{\partial_z\Phi(z)} \quad \text{on } \mathbb{C}.
\end{equation*}
Due to (2.10) we have the global uniform bound $|K_1(z)| < \frac{1}{2}$ everywhere. By combining this with the nonnegativity of the Jacobian determinant, $\det D\Phi \ge 0$ (Lemma 2.2), we obtain the distortion inequality for the Wirtinger derivatives:
\begin{equation*}
|\partial_{\bar{z}}\Phi(z)| \le (|K_1(z)| + |K_2(z)|)|\partial_z\Phi(z)|.
\end{equation*}
We define the local maximal dilatation quotient $D_{\Phi}(z)$. In view of condition (1.5), it satisfies everywhere the upper bound:
\begin{equation*}
D_{\Phi}(z) := \frac{\lvert\partial_z\Phi\rvert + \lvert\partial_{\bar{z}}\Phi\rvert}{\lvert\partial_z\Phi\rvert - \lvert\partial_{\bar{z}}\Phi\rvert} \le \frac{\Lambda}{\lambda} = \frac{g'(|\nabla u|)}{|\nabla u|g''(|\nabla u|)}.
\end{equation*}
Since $g$ is a strictly convex function of linear growth and $u \in C^2(\mathbb{R}^2)$, the maximal dilatation $D_{\Phi}(z)$ is locally bounded. Consequently, $\Phi$ is a classical mapping of bounded distortion. By Reshetnyak's regularity theorem, any non-constant solution of the elliptic system inherits the qualitative topological properties of holomorphic functions. Namely, the gradient mapping \\$\Phi : \mathbb{C} \to \mathbb{C}$ is both open and discrete.

In order to show in more detail that the gradient mapping $\Phi : \mathbb{C} \to \mathbb{C}$ is both open and discrete, we can refer to Theorem 1.1 from the work [11] of
Manfredi and Villamor or to the paper of Iwaniec and Sverak [12]. For quoting  the second reference we just note that a mapping being locally open and discrete clearly inherits these properties globally.\\

Since $\Phi$ is discrete and open, we can apply the generalized topological factorization theorem of Stoilow as for instance it is stated on p. 305 in the work of Luisto and Pankka [13]: there exists a decomposition of the form 
\begin{equation*}
\Phi(z) = \Psi(W(z))
\end{equation*}
where:
\begin{enumerate}
    \item $W : \mathbb{C} \to \Omega \subset \mathbb{C}$ is a global quasiconformal homeomorphism satisfying the Beltrami equation $\partial_{\bar{z}}W(z) = \nu(z)\partial_zW(z)$ everywhere, with the complex dilatation defined by $\nu(z) = \partial_{\bar{z}}\Phi / \partial_z\Phi$ on the set where $\partial_z\Phi \neq 0$.
    \item $\Psi : \Omega \to \mathbb{C}$ is a holomorphic function defined on the image domain $\Omega = W(\mathbb{C})$.
\end{enumerate} 
In addition we refer to the analytic version of the factorization theorem presented in [8], which keeps the degree of regularity.

The maximal dilatation of the global homeomorphism $W$ is pointwise bounded by the dilatation of the gradient, i.e., $D_W(z) \le D_{\Phi}(z)$. For later purposes we recall the Lehto-type asymptotic growth condition: our assumption (1.6) implies the divergence of the integrated dilatation over concentric circles, more precisely we have
\begin{equation}
\int_{1}^{\infty} \frac{1}{rK(r)} \, dr = \infty, \quad \text{where } K(r) := \int_{0}^{2\pi} D_{\Phi}(re^{i\theta}) \, d\theta.
\end{equation}

In order to continue we need to know that the image domain of the quasiconformal homeomorphism W is the complex plane:
\begin{equation}
\Omega = W(\mathbb{C}) = \mathbb{C}.
\end{equation}
However, we can not immediately guarantee that (3.2) holds.
We proceed by contradiction: assume that the image domain $\Omega = W(\mathbb{C})$ is not the entire complex plane, i.e., $\Omega \neq \mathbb{C}$. Since $W : \mathbb{C} \to \Omega$ is a global homeomorphism, $\Omega$ is a proper, simply connected domain in $\mathbb{C}$. By the Riemann mapping theorem (see [6], Theorem 14.8, p. 302) there exists a biholomorphic (conformal) mapping $h : \Omega \to \mathbb{D}$ onto the open unit disk $\mathbb{D} = \{w \in \mathbb{C} : |w| < 1\}$.

Let us define the composite mapping $F : \mathbb{C} \to \mathbb{D}$ as $F(z) = (h \circ W)(z)$. By construction, $F$ is a homeomorphism from $\mathbb{C}$ \textit{onto} $\mathbb{D}$, meaning $F(\mathbb{C}) = \mathbb{D}$. Since $h$ is conformal, the complex dilatation is preserved, and the maximal dilatation satisfies $D_F(z) = D_W(z) \le D_\Phi(z)$ pointwise. Thus, $F$ satisfies the Lehto growth condition (3.1) on $\mathbb{C}$.

For any radius $R > 1$, consider the annulus $A(1, R) = \{z \in \mathbb{C} : 1 < |z| < R\}$. Its image under $F$ is a topological ring domain $\Omega_R = F(A(1, R)) \subset \mathbb{D}$. The complement of $\Omega_R$ in the entire complex plane, $\mathbb{C} \setminus \Omega_R$, consists of exactly two disjoint connected components:
\begin{itemize}
    \item $E_R := F(\overline{\mathbb{D}_1})$, which is the image of the closed unit disk. Since $\overline{\mathbb{D}_1}$ is compact and $F$ is continuous, $E_R$ is a compact, bounded set in $\mathbb{D}$. Evidently, $E_R$ is  independent of $R$. Since $F$ is a homeomorphism, $d = \mathrm{diam}(E_R) > 0$ is a fixed positive constant.
    \item $F_R := (\mathbb{C} \setminus \mathbb{D}) \cup F(\{z \in \mathbb{C} : |z| \ge R\})$. This component is closed and unbounded because it  contains the entire exterior of the unit disk, $\mathbb{C} \setminus \mathbb{D}$.
\end{itemize}

The geometric assumptions of Theorem 20.9.3 are   satisfied for the ring domain $\Omega_R = \mathbb{C} \setminus (E_R \cup F_R)$. Therefore, there exists an increasing homeomorphism $\phi : [0, \infty] \to [0, \infty]$, independent of the continua, such that
\[
\mathrm{Mod}(\Omega_R) \le \phi\left( \frac{\mathrm{dist}(E_R, F_R)}{\mathrm{diam}(E_R)} \right).
\]
Since $E_R \subset \mathbb{D}$ and $\mathbb{C} \setminus \mathbb{D} \subset F_R$, the Euclidean distance between $E_R$ and $F_R$ cannot exceed the distance from $E_R$ to the boundary $\partial\mathbb{D}$.  The unit disk has  diameter  $2$, hence we have $\mathrm{dist}(E_R, F_R) \le 2$. By the monotonicity of $\phi$, we obtain a uniform upper bound independent of $R$:
\begin{equation}\label{eq:geometric_bound}
\mathrm{Mod}(\Omega_R) \le \phi\left( \frac{2}{d} \right) = C < \infty \quad \text{for all } R > 1.
\end{equation}

On the other hand, the  modulus  estimate from below for mappings with integrable dilatation (compare [8], Corollary 20.9.2) provides the following lower bound for the conformal modulus of $\Omega_R$ 
\begin{align}
\mathrm{Mod}(\Omega_R) &= \mathrm{Mod}(F(A(1, R))) \nonumber \\
&\ge \frac{1}{2\pi} \int_1^R \frac{1}{rK(r)} \, dr, \label{eq:analytical_bound}
\end{align}
where the integrated maximal dilatation is defined as
\[
K(r) = \int_0^{2\pi} D_\Phi(re^{i\theta}) \, d\theta.
\]
By the Lehto condition (3.1), the integral on the right-hand side of \eqref{eq:analytical_bound} diverges to infinity as $R \to \infty$. Thus, by taking the limit, we obtain
\begin{equation}\label{eq:analytical_bound_lim}
\lim_{R \to \infty} \mathrm{Mod}(\Omega_R) = \infty.
\end{equation}

The divergence stated in \eqref{eq:analytical_bound_lim} contradicts the uniform upper bound established in \eqref{eq:geometric_bound}, which requires $\mathrm{Mod}(\Omega_R) \le C$ to hold globally for all $R > 1$.

 Consequently, the initial assumption $\Omega \neq \mathbb{C}$ was false, which completes the proof of
\[
\Omega = W(\mathbb{C}) = \mathbb{C}.
\]

Now, on account of (3.2), $\Psi(z) = \Phi(W^{-1}(z))$ is defined for any complex number z and  thereby a classical entire holomorphic function. In consequence, Picard's little theorem applies to $\Psi$ and we arrive at the "or" case from Theorem 1.1.
\qed
\section{Remarks concerning Theorem 1.6}
The coefficient matrix in this new setting is just the Hessian matrix of $f$: we define
\begin{equation*}
a_{ij}(z) := \frac{\partial^{2} f}{\partial P_{i}\partial P_{j}}(\nabla u(z)) \quad (i,j = 1,2).
\end{equation*}
Inequality (2.5) requires the strict positivity of the trace of the coefficient matrix, i.e., $a_{11}+a_{22}>0$. Since $a_{ij}$ represents the components of the Hessian matrix $D^{2}f(\nabla u)$, the sum $a_{11}+a_{22}$ corresponds exactly to $\text{trace}(D^{2}f)$. From condition \eqref{eq:2}, it follows for the eigenvalues $\lambda_{1}, \lambda_{2}$ of $D^{2}f$ that:
\begin{equation*}\label{eq:4}
a_{11}+a_{22} = \lambda_{1} + \lambda_{2} \ge 2c_{1}(1+|\nabla u|)^{-\mu} > 0
\end{equation*}
Since $c_{1} > 0$, this sum is strictly positive everywhere. Consequently, (2.5) holds in the new setting.\\
The transformation of the Euler--Lagrange equation (1.8) into the generalized Beltrami equation
\begin{equation*}
\partial_{\overline{z}}\Phi = K_{1}\partial_{z}\Phi + K_{2}\overline{\partial_{z}\Phi}
\end{equation*}
is of purely algebraic nature: expanding equation (1.8) $\text{div}\,[Df(\nabla u)] = 0$ via the chain rule yields, due to the symmetry of the mixed derivatives ($f_{P_{1}P_{2}} = f_{P_{2}P_{1}}$):
\begin{equation*}\label{eq:5}
a_{11}p_{x} + 2a_{12}p_{y} + a_{22}q_{y} = 0
\end{equation*}
This matches equation (2.3) exactly. The proof of Lemma 2.1 thus remains valid without modification, and the coefficients $K_{1}, K_{2}$ possess the identical structure as in (2.8).
We check the validity of Lemma 2.2 and inequality (2.10):
\begin{itemize}
    \item \textbf{Lemma 2.2 ($\det D\Phi \ge 0$):} The proof just uses relation (2.3) and the positive definiteness of the coefficient matrix. Since condition \eqref{eq:2} guarantees the positive definiteness of $D^{2}f$, the Jacobian determinant remains non-negative.
    
    \item \textbf{Equation (2.10) ($|K_{1}| < \frac{1}{2}$):} Due to the real-valuedness of $f$, the coefficients $a_{ij}$ are real, which implies $K_{2} = \overline{{K}_{1}}$. For the absolute value, it holds that:
    \begin{equation*}\label{eq:6}
    |K_{1}|^{2} = \frac{1}{4} - \frac{\det(a_{ij})}{(\text{trace}(a_{ij}))^{2}} = \frac{1}{4} - \frac{\lambda_{1}\lambda_{2}}{(\lambda_{1} + \lambda_{2})^{2}}
    \end{equation*}
    By the strict ellipticity, the eigenvalues satisfy $\lambda_{1}, \lambda_{2} > 0$, and consequently $\lambda_{1}\lambda_{2} > 0$. Thus, a strictly positive value is subtracted from $\frac{1}{4}$, ensuring that $|K_{1}| < \frac{1}{2}$ is always satisfied.
\end{itemize}
Since $\det D\Phi \ge 0$ and $|K_{1}| < \frac{1}{2}$, $\Phi$ is a mapping of bounded distortion, and as before it can be factored as $\Phi(z) = \Psi(W(z))$ with a quasiconformal homeomorphism $W$ and a holomorphic function $\Psi$.

The maximal dilatation satisfies:
\begin{equation*}\label{eq:7}
D_{\Phi}(z) \le \frac{\lambda_{\max}(D^{2}f)}{\lambda_{\min}(D^{2}f)} \le \frac{c_{2}(1+|\nabla u|)^{-1}}{c_{1}(1+|\nabla u|)^{-\mu}} = \frac{c_{2}}{c_{1}}(1+|\nabla u|)^{\mu-1}
\end{equation*}
This corresponds precisely to the asymptotic growth behavior found in the radial case. Under the assumption of the growth condition (1.6) for the gradient, the integrated dilatation grows logarithmically, whereby the Lehto condition is satisfied and $\Omega = W(\mathbb{C}) = \mathbb{C}$ follows, provided we assume that $\Phi$ is not constant. Thus, the arguments in Section 3 carry over completely.
\qed
\section{Proof of Theorem 1.8}
As in Section 3 we argue by contradiction assuming that the complex gradient mapping $\Phi$ is not constant.  With the help of (1.10) we  will derive a precise asymptotic lower bound for the global quasiconformal homeomorphism $W(z)$ for large values of $|z|$ leading to the desired contradiction. 

From the former growth condition (1.6) for the gradient,
\[
|\nabla u(z)| \le c \, [\ln(2 + |z|)]^{\frac{1}{\mu-1}}, \quad z \in \mathbb{R}^2,
\]
combined with the previous dilatation estimate, we deduce that the maximal dilatation quotient satisfies an asymptotic upper bound of the form $D_\Phi(z) \le C \ln|z|$ for sufficiently large $|z|$. Consequently, the integrated maximal dilatation $K(r)$ defined in (3.1) grows at most logarithmically, more precisely it holds, 
\[
K(r) = \int_0^{2\pi} D_\Phi(re^{i\theta}) \, d\theta \le C_0 \ln r \quad \text{as } r \to \infty,
\]
where $C_0 > 0$ is a positive constant depending only on the data.

Using the distortion property under the global quasiconformal mapping $W$, the conformal modulus of the image of the annulus $A(1, R) = \{z \in \mathbb{C} : 1 < |z| < R\}$ satisfies the  lower Lehto estimate (compare again [8], Corollary 20.9.2)
\begin{equation}\label{eq:lehto_lower}
\mathrm{Mod}(W(A(1, R))) \ge \frac{1}{2\pi} \int_1^R \frac{1}{rK(r)} \, dr.
\end{equation}
Substituting the logarithmic upper bound for $K(r)$ into \eqref{eq:lehto_lower} yields a double-logarithmic lower estimate for the modulus of the form
\begin{align}
\mathrm{Mod}(W(A(1, R))) &\ge \frac{1}{2\pi} \int_1^R \frac{1}{C_0 \, r \ln r} \, dr \nonumber \\
&= \frac{1}{2\pi C_0} \ln(\ln R). \label{eq:modulus_double_log}
\end{align}

On the other hand, we can normalize \(W(0)=0\) and choose a radius \(\rho>0\) such that
\(\overline{B(0,\rho)}\subset W(B(0,1))\). We further abbreviate \(m(R)=\min_{|z|=R}|W(z)|\). The Grötzsch ring estimate, applied after inversion, gives
\begin{equation}\label{equ:upper bound}
\operatorname{Mod}(W(A(1,R))) \leq \frac{1}{2\pi}\ln\frac{4m(R)}{\rho}.
\end{equation}
For details the reader is referred to the monograph [14] of Gehring, Martin and Palka.

Combining the lower bound \eqref{eq:modulus_double_log} with the upper bound \eqref{equ:upper bound} leads to
\[\frac{1}{2\pi C_0} \ln(\ln R) \le \frac{1}{2\pi} \ln \left( \min_{|z|=R} |W(z)| \right),\] where for notational simplicity we have just assumed that \(\rho=4.\)
 Applying the exponential function to both sides, the following lower bound for $|z| = R$ sufficiently large is established:
\begin{equation}\label{eq:repaired_lower_bound}
|W(z)| \ge c_1 \, (\ln |z|)^\alpha \quad \text{for } |z| \gg 1,
\end{equation}
where $\alpha = \frac{1}{C_0} > 0$ and $c_1 > 0$ are positive constants. 

Inequality \eqref{eq:repaired_lower_bound} immediately translates into an upper bound for the inverse mapping $W^{-1}$, more precisely, it holds
\begin{equation}\label{eq:repaired_inverse_bound}
|W^{-1}(z)| \le \exp \left( \left( \frac{|z|}{c_1} \right)^{\frac{1}{\alpha}} \right) \quad \text{for } |z| \gg 1.
\end{equation}

Recalling the Stoilow factorization $\Phi(z) = \Psi(W(z))$, the entire holomorphic function $\Psi : \mathbb{C} \to \mathbb{C}$ is given by $\Psi = \Phi \circ W^{-1}$. Combining the gradient growth condition (1.10) with the  upper bound \eqref{eq:repaired_inverse_bound}, we obtain an estimate for the modulus of $\Psi(z)$:

\begin{equation}
\begin{aligned}
\lvert\Psi(z)\rvert &\le c \ln \left[ \ln \left\lvert W^{-1}(z) \right\rvert \right] 
\le c \ln \left( \left[ \frac{\lvert z\rvert}{c_1} \right]^{\frac{1}{\alpha}} \right) \\
&= \frac{c}{\alpha} \ln \left( \frac{\lvert z\rvert}{c_1} \right) \le c_2 \ln \lvert z\rvert
\end{aligned}
\end{equation}
being valid for sufficiently large values of $z$.
 By a variant of Liouville's theorem we deduce from (5.6) the constancy of $\Psi$ and thereby of $\Phi$, leading to the desired contradiction, which  completes the proof of Theorem 1.8.
\qed
\section{Proof of Theorem 1.10}

Assuming that the conditions (1.2)--(1.5) hold for the energy density $g$ with a fixed exponent $\mu > 1$, we consider an entire solution $u \in C^2(\mathbb{R}^2)$ to the Euler--Lagrange equation (1.1) or its more general variant (1.8) under the $\mu$-ellipticity condition (1.9). 

Instead of the  Lehto growth condition (1.6) or the double-logarithmic bound (1.10) required in Theorem 1.8, we now impose the Landau $o$-condition (1.11) for the gradient at infinity, i.e., we require
\begin{equation}
|\nabla u(z)| = o\left([\ln(2 + |z|)]^{\frac{1}{\mu-1}}\right) \quad \text{as } |z| \to \infty,
\end{equation}
and claim that u must be an affine function.

By contradiction, suppose that the complex gradient mapping $\Phi(z) = u_x - iu_y$ is non-constant. As outlined in Section 3 and 5 $\Phi$ is an open and discrete mapping, and using Stoilow's factorization theorem, $\Phi$ decomposes into $\Phi = \Psi \circ W$, where $W: \mathbb{C} \to \mathbb{C}$ is a global quasiconformal homeomorphism and $\Psi: \mathbb{C} \to \mathbb{C}$ is an entire holomorphic function. We now repeat the calculations from Section 5 having a more careful look at the constants:

Under the $o$-condition (6.1), the local maximal dilatation satisfies $D_\Phi(z) \le \varepsilon C \ln|z|$ for any choice of $\varepsilon > 0$ and sufficiently large $|z|$. Consequently, the integrated dilatation over concentric circles satisfies the asymptotic upper bound:
\begin{equation}
K(r) = \int_0^{2\pi} D_\Phi(re^{i\theta}) \, d\theta \le \varepsilon C_0 \ln r \quad \text{as } r \to \infty.
\end{equation}
Inserting  estimate (6.2) into the lower Lehto distortion inequality for the conformal modulus of the ring domain yields:
\begin{equation}
\text{Mod}(W(A(1, R))) \ge \frac{1}{2\pi \varepsilon C_0} \ln(\ln R).
\end{equation}
Combining (6.3) with the upper bound derived from the Grötzsch ring estimate (compare (5.3) and again choose \(\rho=4\)) , we apply the exponential function to both sides. This implies a controlled asymptotic lower bound for the global homeomorphism $W(z)$ for $|z| \gg 1$:
\begin{equation}
|W(z)| \ge c_1(\ln|z|)^{\alpha_\varepsilon} \quad \text{with } \alpha_\varepsilon = \frac{1}{\varepsilon C_0}.
\end{equation}
Inequality (6.4) yields a sub-exponential growth rate for the inverse mapping $W^{-1}$ for $|z| \gg 1$, more precisely, it holds
\begin{equation}
|W^{-1}(z)| \le \exp\left( \left(\frac{|z|}{c_1}\right)^{\frac{1}{\alpha_\varepsilon}} \right) = \exp\left( \left(\frac{|z|}{c_1}\right)^{\varepsilon C_0} \right).
\end{equation}
We now estimate the modulus of the entire holomorphic function $\Psi(z) = \Phi(W^{-1}(z))$ by inserting the growth bound (6.5)  into inequality (1.6):
\begin{equation}
|\Psi(z)| \le c \left[ \ln\left(2 + \exp\left( \left(\frac{|z|}{c_1}\right)^{\varepsilon C_0} \right)\right) \right]^{\frac{1}{\mu-1}} \le C_3 |z|^{\frac{\varepsilon C_0}{\mu-1}} = C_3 |z|^\gamma.
\end{equation}
The resulting growth exponent is defined as $\gamma := \frac{\varepsilon C_0}{\mu - 1}$. Since $\varepsilon > 0$ can be chosen arbitrarily small from the very beginning, the requirement for a strictly sub-polynomial growth rate\begin{equation*}
\gamma < 1
\end{equation*}
can evidently be fulfilled. By the extended variant of Liouville's theorem for entire functions, the sub-polynomial growth  stated in (6.6)  shows that $\Psi$, and thereby $\Phi$, is a constant vector. This provides the desired contradiction to our initial assumption, completing the proof that $u$ must be an affine function.\qed
\begin{remark}
We think that the proof of Theorem 1.10 can be adjusted to the case of densities $f$  satisfying an anisotropic  $(p,q)$-ellipticity condition with exponents $2 \le p \le q$ in place of (1.9). If $p<q$, we then have to replace (1.11) by the condition 
  \begin{equation*}
    |\nabla u(z)| = o\left([\ln(2 + |z|)]^{\frac{1}{q-p}}\right) \quad \text{for } |z| \to \infty
    \end{equation*} 
    implying that $u$ is an affine function, whereas in the isotropic setting $q=p$ this follows from the requirement $|\nabla u(z)| = \mathcal{O}(|z|^\alpha)$ with the sharp bound for the exponent
    \begin{equation*}
    \alpha < \frac{1}{p-1}.
    \end{equation*}
\end{remark}
\section*{Acknowledgements}The author is very grateful to Kari Astala and Gaven Martin for their helpful comments.

\section*{}Funding. Open Access funding enabled and organized by Projekt DEAL.\\\\
Data availability. No additional data are involved.\\\\
Declarations\\
Conflict of interest.  There is no
conflict of interest to be declared.\\\\\\
\begin{minipage}[t]{0.48\textwidth}
\textbf{Martin Fuchs}\\
Saarland University\\
Department of Mathematics\\
P.O. Box 15 11 50\\
66041 Saarbrücken, Germany\\
\texttt{fuchs@math.uni-sb.de}
\end{minipage}
\end{document}